\documentclass[review]{elsarticle}
\usepackage{amsmath, amssymb, amsthm}
\usepackage{textcomp}
\newtheorem{theorem}{Theorem}[section]
\newtheorem{lemma}{Lemma}[section]
\newtheorem{definition}{Definition}[section]
\newtheorem*{remark}{Remark}

\journal{Journal of \LaTeX\ Templates}

\begin{document}

\begin{frontmatter}

\title{ Compactness of Semicommutators of Toeplitz operators- a Characterization}

%\author{V.B. Kiran Kumar}
%\ead{kiranbalu36@gmail.com}
%\author{M.N.N. Namboodiri}
%\ead{mnnadri@gmail.com}
\author{Rahul Rajan}
\ead{rajan.rahul48@yahoo.com}
\address{Department of Mathematics, Cochin University of Science and Technology, Cochin, India}


\begin{abstract}
 Let $T_{f}$ denote the Toeplitz operator on the Hardy space $H^{2}(\mathbb{T})$ and let $T_{n}(f)$ be the corresponding $n \times n$ Toeplitz matrix. In this paper, we characterize the compactness of the operators $T_{|f|^{2}}-T_{f}T_{\overline{f}}$ and $T_{|\tilde{f}|^{2}}-T_{\tilde{f}}T_{\overline{\tilde{f}}},$ where $\tilde{f}(z)=f(z^{-1}),$ in terms of the convergence of the sequence $\{T_{n}(|f|^{2})-T_{n}(f)T_{n}(\overline{f})\}$ in the sense of singular value clustering. Hence we obtain a method to check the compactness of semicommutators of Toeplitz operators using the matrices obtained from the Fourier coefficients of the symbol function (Toeplitz matrices). The function space $VMO \cap L^{\infty}(\mathbb{T})$ is the largest $C^{*}$-subalgebra of $L^{\infty}(\mathbb{T})$ with the property that whenever $f,g \in VMO \cap L^{\infty}(\mathbb{T})$, $T_{fg}-T_{f}T_{g}$ is compact. In this article, we obtain a characterization of $VMO \cap L^{\infty}(\mathbb{T})$ in terms of the convergence of $\{T_{n}(fg)-T_{n}(f)T_{n}(g)\}$ in the sense of singular value clustering. To be precise, $VMO \cap L^{\infty}(\mathbb{T})$ is the largest $C^{*}$-subalgebra of $L^{\infty}(\mathbb{T})$ with the property that whenever $f,g \in VMO \cap L^{\infty}(\mathbb{T})$, $\{T_{n}(fg)-T_{n}(f)T_{n}(g)\}$ converges in the sense of singular value clustering.
\end{abstract}

\begin{keyword}
Toeplitz Operators \sep Toeplitz Matrices
%\MSC[2010] 00-01\sep  99-00
\end{keyword}

\end{frontmatter}

%\linenumbers
\section{Introduction}
Let $\mathbb{T}$ be the unit circle. Consider the subspace $H^{2}(\mathbb{T})$ of $L^{2}(\mathbb{T}, \frac{1}{2\pi}d\theta)$ defined as follows.
$$ H^{2}(\mathbb{T})=\{ f \in L^{2}(\mathbb{T}, \frac{1}{2\pi}d\theta)\,\, :\,\, f_{n}=0\textrm{ for all $n <0$} \}, $$
where $f_{n}$ is the $n^{th}$ Fourier coefficient of $f$. Let $P$ be the projection of $L^{2}(\mathbb{T}, \frac{1}{2\pi}d\theta)$ onto $H^{2}(\mathbb{T})$. For $f \in L^{2}(\mathbb{T}, \frac{1}{2\pi}d\theta),$ the Toeplitz operator $T_{f}$ defined on $H^{2}(\mathbb{T})$ is as follows.
$$ T_{f}(g)=P(fg),\,\,g\in H^{2}(\mathbb{T}). $$
The domain of $T_{f}$ contains $H^{\infty}$, the set of all essentially bounded functions which are analytic on $\mathbb{T}$. Note that $H^{\infty}$ is dense in $H^{2}(\mathbb{T})$. It is known that $T_{f}$ is bounded if and only if $f$ is bounded. With respect to the the orthonormal basis $\{z^{n}\,\, :\,\, n \geq 0\}$, the bounded Toeplitz operator has the matrix representation
$$ (T_{f})_{i,j}=(f_{i-j}),\,\, 0\leq i,j < \infty. $$ 
For $f \in L^{2}(\mathbb{T}, \frac{1}{2\pi}d\theta)$, we shall identify the operator $T_{f}$ with the infinite matrix $(f_{i-j})$.
The Toeplitz matrix $T_{n}(f)$ is the matrix $(f_{i-j}),\,\, 0\leq i,j < n.$ Also, it is easy to observe that $T_{n}(f)=P_{n}T_{f}P_{n}$, where $P_{n}$ is the projection defined by $P_{n}(x_{0},x_{1},\hdots,x_{n-1},x_{n},\hdots)=(x_{0},x_{1},\hdots,x_{n-1},0,0,\hdots)$.\\

Toeplitz matrices and operators arise in and have application to a wide variety
of fields of pure and applied mathematics, particularly, numerical analysis,
probability theory, complex analysis, K-theory, statistical mechanics, etc.\\

For $f \in L^{2}(\mathbb{T},\frac{1}{2\pi}d\theta)$, the Hankel operator $H_{f}$ is defined as follows.
$$ H_{f}: H^{2}(\mathbb{T})\rightarrow (H^{2}(\mathbb{T}))^{\perp}, $$
$$ H_{f}(g)=(I-P)(fg). $$
Note that 
$$ T_{fg}-T_{f}T_{g}=H_{\bar{f}}^{*}H_{g}. $$
If $f$ is bounded, then clearly, $H_{f}$ is bounded with $\|H_{f}\|\leq \|f\|_{\infty}$. In general, $H_{f}$ is densely defined with its domain containing $H^{\infty}$. There are unbounded symbols $f$ for which $H_{f}$ is bounded. Consider the class of functions $BMO$ and $VMO$ which is defined as follows ($BMO$ stands for Bounded Mean Oscillation and $VMO$ stands for Vanishing Mean Oscillation).
 
\begin{definition}\cite{sarson}
 Let $f\in L^{2}(\mathbb{T},\frac{1}{2\pi}d\theta)$. Consider an interval $I$ on $\mathbb{T}$. Define,
 $$ f_{I}=\frac{1}{|I|} \int_{I}f(\theta)d\theta, $$
where $|I|$ is the length of the interval $I$. $f$ is said to have bounded mean oscillation or $f\in BMO$ if,
 $$ \sup_{I} \bigg[\frac{1}{|I|}\int_{I}|f(\theta)-f_{I}|^{2}d\theta \bigg]^{1/2} < \infty. $$\\
 
 $f$ is said to have vanishing mean oscillation or $f\in VMO$ if,
 $$ \lim_{|I|\rightarrow0} \bigg[\frac{1}{|I|}\int_{I}|f(\theta)-f_{I}|^{2}d\theta \bigg]^{1/2}=0. $$\\
\end{definition}
$BMO$ characterizes the functions $f$ such that the Hankel operators $H_{f}$ and $H_{\bar{f}}$ are simultaneously bounded. Also, there are unbounded symbols $f$ for which $H_{f}$ is compact. $VMO$ characterizes the functions $f$ such that both $H_{f}$ and $H_{\bar{f}}$ are compact. 
\begin{remark}
Note that $L^{\infty}(\mathbb{T})\subsetneq BMO$ and $C(\mathbb{T})\subsetneq VMO$. Also,
 $VMO \cap L^{\infty}(\mathbb{T})$ is a $C^{*}$-subalgebra of $L^{\infty}(\mathbb{T})$ such that whenever $f,g \in VMO \cap L^{\infty}(\mathbb{T}),$ $T_{fg}-T_{f}T_{g}$ is compact. In fact, this is the largest $C^{*}$-subalgebra with this property (see \cite{zhu}).
\end{remark}

The problem of checking the compactness of $T_{|f|^{2}}-T_{f}T_{\bar{f}}$ is one of the important questions in operator theory. There are several characterizations of this in terms of symbols (see \cite{zhu}). Here we obtain a characterization using the matrices formed by the Fourier coefficients of the symbol. The advantage is that we require only the information of Fourier coefficients of the symbols to check the compactness of the semicommutators. Now we define the convergence of a sequence of matrices $\{A_{n}\}$, where each $A_{n}$ is a $n \times n$ matrix, in the sense of singular value clustering.

\begin{definition}\cite{laa}
 Let $\{A_{n}\}$ be a sequence of matrices where each $A_{n}$ is a matrix of order $n \times n$. We say that $\{A_{n}\}$ converges to $\{O_{n}\}$ in Type 2 strong cluster sense, if for $\epsilon >0$, there exists two positive integers $N_{1,\epsilon}=O(1)$ and $N_{2,\epsilon}$ such that for $n>N_{2,\epsilon}$, 
 $$ A_{n}=R_{n}+N_{n}, $$
 where $Rank(R_{n}) \leq N_{1,\epsilon}$ and $\|N_{n}\|< \epsilon$. Also, if $N_{1,\epsilon}=o(n)$, then we say $\{A_{n}-B_{n}\}$ converges to $\{O_{n}\}$ in Type 2 weak cluster sense.
  \end{definition}
 \begin{remark}
 There is a notion of the convergence in Type 1 strong cluster sense (Type 1 weak cluster sense respectively). $\{A_{n}\}$ converges to $\{O_{n}\}$ in Type 1 strong cluster sense if for each $\epsilon>0$, there exists 2 positive integers $N_{1,\epsilon}=O(1)$ ($o(n)$ respectively) and $N_{2,\epsilon}$ such that for $n>N_{2,\epsilon}$, except atmost $N_{1,\epsilon}$ eigenvalues, all other eigenvalues of $A_{n}$ belong to the $\epsilon$ neighborhood of $0$. Note that if $A_{n}$ is normal for each $n$, then Type 1 and Type 2 convergence are equivalent (see \cite{laa}).\\
 
 The notion of convergence in Type1/Type2 strong cluster sense arises from the preconditioning problem in numerical linear algebra (see \cite{ser,tr}). The convergence in Type 2 strong cluster sense is equivalent to the singular value clustering as we see in the following result. 
\end{remark}
\begin{lemma}\cite{laa}
 $\{A_{n}\}$ converges to $\{O_{n}\}$ in Type 2 strong cluster sense if and only if for each $\epsilon>0$, there exists 2 positive integers $N_{1,\epsilon}=O(1)$ and $N_{2,\epsilon}$ such that for $n>N_{2,\epsilon}$, except atmost $N_{1,\epsilon}$ singular values, all other singular values of $A_{n}$ belong to the interval $[0,\epsilon)$.
\end{lemma}
\subsection{Main Results}
In this article, we prove the following results.
\begin{theorem}\label{eq}
 For $f \in BMO$, the following are equivalent.
 \begin{itemize}
  \item[1.] $\{T_{n}(|f|^{2})-T_{n}(f)T_{n}(\bar{f})\}$ converges to $\{O_{n}\}$ in Type 2 strong cluster sense.
  \item[2.] $T_{|f|^{2}}-T_{f}T_{\bar{f}}$ and $T_{|\tilde{f}|^{2}}-T_{\tilde{f}}T_{\tilde{\bar{f}}}$ are compact 
 \end{itemize}

\end{theorem}
\begin{theorem}\label{vmo}
 Let $f\in BMO$. Then $f \in VMO$ if and only if $\{T_{n}(|f|^{2})-T_{n}(f)T_{n}(\bar{f})\}$ and $\{T_{n}(|f|^{2})-T_{n}(\bar{f})T_{n}(f)\}$ converge to $\{O_{n}\}$ in Type 2 strong cluster sense. 
\end{theorem} 
 \begin{theorem}\label{vmoc}
 Whenever $f,g \in VMO \cap L^{\infty}(\mathbb{T})$, $\{T_{n}(fg)-T_{n}(f)T_{n}(g)\}$ converges to $\{O_{n}\}$ in Type 2 strong cluster sense. Also, $VMO \cap L^{\infty}(\mathbb{T})$ is the largest $C^{*}$-subalgebra of $L^{\infty}(\mathbb{T})$ with this property. 
\end{theorem}

 The advantage of these results is that we use the Fourier coefficients of $f$ to derive the properties of operators formed by $f$. This is in line with the classical Fourier analysis where one tries to retrieve the symbol $f$ from its Fourier coefficients. Note that, one shall form a Toeplitz matrix by knowing the Fourier coefficients of its symbol.  It was obtained that an operator $K \in \mathcal{B}(\mathcal{H})$, the space of all bounded linear operators on a separable Hilbert space $\mathcal{H}$, is compact if and only if $\{K_{n}=P_{n}KP_{n}\}$ converges to $\{O_{n}\}$ in Type 2 strong cluster sense. Here, we shall check the convergence of $\{T_{n}(|f|^{2})-T_{n}(f)T_{n}(\bar{f})\}$ instead of $\{P_{n}(T_{|f|^{2}}-T_{f}T_{\bar{f}})P_{n}\}$. The former one is comparatively easy to compute, as we need to know the Fourier coefficients of $f$, to determine the matrix $T_{n}(|f|^{2})-T_{n}(f)T_{n}(\bar{f})$.\\
 
 Theorem \ref{eq} establishes a connection between  compactness, an operator theoretic notion, and the convergence of matrices in strong cluster sense, a notion in numerical linear algebra. Theorem \ref{vmo} and Theorem \ref{vmoc} establish a connection between the function space $VMO \cap L^{\infty}(\mathbb{T})$ and the convergence in strong cluster sense. In short, these results establish a connection between various notions in operator theory, complex function spaces, and numerical linear algebra.\\
 
\section{Proof of Main Results}

First, we give a few results which are useful to prove our main results.
\begin{lemma}\label{wid}\cite{wid}
 $$ T_{n}(fg)-T_{n}(f)T_{n}(g)=P_{n}H_{\bar{f}}^{*}H_{g}P_{n}+Q_{n}H_{\tilde{\bar{f}}}^{*}H_{\tilde{g}}Q_{n}, $$
 where $Q_{n}$ is the operator defined by 
 $$ Q_{n}(x_{0},x_{1},\hdots,x_{n-1},x_{n},\hdots)=(x_{n-1},x_{n-2},\hdots,x_{0},0,0,\hdots) $$ and $\tilde{f}(z)=f(z^{-1})$.
\end{lemma}
\begin{remark}
 In \cite{bott,wid}, the Hankel operator is defined as follows.
 $$ H(a) : H^{2}(\mathbb{T})\rightarrow H^{2}(\mathbb{T}), $$
 $$ f \mapsto PM(a)(I-P)J(f), $$ 
 where $M(a)$ denotes the Multiplication operator corresponding to $a$ and $J(f)(z)=z^{-1}f(z^{-1})$. Also, note that the matrix representation of $H(a)$ and $H_{\bar{a}}^{*}$ are the same. In \cite{wid}, Lemma \ref{wid} is equivalently stated as follows.
 $$ T_{n}(fg)-T_{n}(f)T_{n}(g)=P_{n}H(f)H(\tilde{g})P_{n}+Q_{n}H(\tilde{f})H(g)Q_{n}. $$
 
\end{remark}
\begin{lemma}\label{pos}\cite{laa}
 Let $\{A_{n}\}$ and $\{B_{n}\}$ be sequences of positive matrices, where each $A_{n}$ and $B_{n}$ are $n \times n$ matrices. If $\{A_{n}+B_{n}\}$ converges to $\{O_{n}\}$ in Type 1 strong cluster sense, then $\{A_{n}\}$ and $\{B_{n}\}$ converge to $\{O_{n}\}$ in Type 1 strong cluster sense. 
\end{lemma}

\begin{lemma}\label{vbk2}\cite{mon,vbk}
 Let $A \in \mathcal{B}(\mathcal{H})$. $A$ is compact if and only if $\{A_{n}\}$ converges to $\{O_{n}\}$ in Type 2 strong cluster sense.
\end{lemma}
\begin{lemma}\label{zhu}\cite{zhu}
 Let $f\in L^{2}(\mathbb{T},\frac{1}{2\pi},d\theta)$. $H_{f}$ and $H_{\bar{f}}$ is bounded if and only if $f\in BMO$. Also, $H_{f}$ and $H_{\bar{f}}$ is compact if and only if $f\in VMO$.  
\end{lemma}

Now we give the proofs of our main results.\\
\textbf{Proof of Theorem \ref{eq}:-}
\begin{proof}
  By Lemma \ref{wid},
 \begin{equation}\label{1}
  T_{n}(|f|^{2})-T_{n}(f)T_{n}(\bar{f})=P_{n}H_{\bar{f}}^{*}H_{\bar{f}}P_{n}+Q_{n}H_{\tilde{\bar{f}}}^{*}H_{\tilde{\bar{f}}}Q_{n}.
  \end{equation}
  Note that $T_{n}(|f|^{2})-T_{n}(f)T_{n}(\bar{f}),\, P_{n}H_{\bar{f}}^{*}H_{\bar{f}}P_{n}$ and $Q_{n}H_{\tilde{\bar{f}}}^{*}H_{\tilde{\bar{f}}}Q_{n}$ are positive matrices.
 Hence by Lemma \ref{pos}, $\{T_{n}(|f|^{2})-T_{n}(f)T_{n}(\bar{f})\}$ converges to $\{O_{n}\}$ in Type 2 strong cluster sense if and only if
 $\{P_{n}H_{\bar{f}}^{*}H_{\bar{f}}P_{n}\}$ and $\{Q_{n}H_{\tilde{\bar{f}}}^{*}H_{\tilde{\bar{f}}}Q_{n}\}$ converge to $\{O_{n}\}$ in Type 2 strong cluster sense. Note that $Q_{n}H_{\tilde{\bar{f}}}^{*}H_{\tilde{\bar{f}}}Q_{n}$ and $P_{n}H_{\tilde{\bar{f}}}^{*}H_{\tilde{\bar{f}}}P_{n}$ are unitarily equivalent and thus $\{P_{n}H_{\tilde{\bar{f}}}^{*}H_{\tilde{\bar{f}}}P_{n}\}$ converges to $\{O_{n}\}$ in Type 2 strong cluster sense if and only if $\{Q_{n}H_{\tilde{\bar{f}}}^{*}H_{\tilde{\bar{f}}}Q_{n}\}$ converges to $\{O_{n}\}$ in Type 2 strong cluster sense. By Lemma \ref{vbk2}, $\{P_{n}H_{\bar{f}}^{*}H_{\bar{f}}P_{n}\}$ converges to $\{O_{n}\}$ if and only if $H_{\bar{f}}^{*}H_{\bar{f}}$ is compact. Similarly, $\{P_{n}H_{\tilde{\bar{f}}}^{*}H_{\tilde{\bar{f}}}P_{n}\}$ converges to $\{O_{n}\}$ if and only if $H_{\tilde{\bar{f}}}^{*}H_{\tilde{\bar{f}}}$ is compact. We know that $T_{|f|^{2}}-T_{f}T_{\bar{f}}=H_{\bar{f}}^{*}H_{\bar{f}}$ and $T_{|\tilde{f}|^{2}}-T_{\tilde{f}}T_{\tilde{\bar{f}}}=H_{\tilde{\bar{f}}}^{*}H_{\tilde{\bar{f}}}$ and hence we obtain the result. 
\end{proof}
\begin{remark}\label{rem1}
 If $f \in VMO$, then $\tilde{f}$ and $\bar{f}$ belong to $VMO$. Hence for $f\in VMO$, the compactness of $T_{|f|^{2}}-T_{f}T_{\bar{f}}$ and the convergence of $\{T_{n}(|f|^{2})-T_{n}(f)T_{n}(\bar{f})\}$ to $\{O_{n}\}$ in Type 2 strong cluster sense are equivalent.
\end{remark}

\textbf{Proof of Theorem \ref{vmo}:- }
\begin{proof}
 Let $f\in VMO$. Then it is easy to check that $\tilde{f}$ and $\bar{f}$ belong to $VMO$. Hence $H_{f}$, $H_{\tilde{f}}$, $H_{\tilde{\bar{f}}}$ and $H_{\bar{f}}$ are compact. Consider (\ref{1}) and
 \begin{equation}\label{2}
  T_{n}(|f|^{2})-T_{n}(\bar{f})T_{n}(f)=P_{n}H_{f}^{*}H_{f}P_{n}+Q_{n}H_{\tilde{f}}^{*}H_{\tilde{f}}Q_{n}.
 \end{equation}
 By Lemma \ref{vbk2}, $\{T_{n}(|f|^{2})-T_{n}(f)T_{n}(\bar{f})\}$ and $\{T_{n}(|f|^{2})-T_{n}(\bar{f})T_{n}(f)\}$ converge to $\{O_{n}\}$ in Type 2 strong cluster sense.\\
 
 Conversely, assume $\{T_{n}(|f|^{2})-T_{n}(f)T_{n}(\bar{f})\}$ and $\{T_{n}(|f|^{2})-T_{n}(\bar{f})T_{n}(f)\}$ converge to $\{O_{n}\}$ in Type 2 strong cluster sense. From (\ref{1}) and (\ref{2}), using Lemma \ref{pos}, we have $\{P_{n}H_{\bar{f}}^{*}H_{\bar{f}}P_{n}\}$ and $\{P_{n}H_{f}^{*}H_{f}P_{n}\}$ converge to $\{O_{n}\}$ in Type 2 strong cluster sense. Thus by Lemma \ref{vbk2}, $H_{f}$ and $H_{\bar{f}}$ are compact and hence $f \in VMO.$ 
\end{proof}
Now we give Uchiyama's inequality which is useful to prove Theorem \ref{vmoc}.\\

A \textit{Schwarz map} on a $C^{*}$-algebra $\mathbb{A}$ is a positive linear map $\psi $ satisfying the condition $\psi(a^{*}a) \geq \psi(a^{*})\psi(a)$, for all $a \in \mathbb{A}$. A \textit{generalized Schwarz map} on $\mathbb{A}$ 
with respect to the binary operation $\circ$ is a map $\Phi$ which satisfies
$\Phi(x^*) = \Phi(x)^*$ and $\Phi(x^*) \circ \Phi(x) \leq \Phi(x^* \circ x)$ for every $x \in \mathbb{A}$. Here '$\circ$' denotes a binary operation on $\mathbb{A}$ with certain properties (see \cite{uch}). We avoid listing these properties, instead we mention that these properties are satisfied by composition of operators, Jordan product and pointwise product of functions. 
\begin{lemma}\label{uch}\cite{uch}
Let $\varPhi$ be a generalized Schwarz map with respect to a binary operation $\circ$ on a $C^{*}$-algebra $\mathbb{A}$. For $f,g \in \mathbb{A}$, let
 $$ X=\varPhi(f^{*}\circ f)-\varPhi(f)^{*} \circ \varPhi(f) \geq 0, $$
 $$ Y=\varPhi(g^{*}\circ g)-\varPhi(g)^{*} \circ \varPhi(g) \geq 0, $$
 $$ Z=\varPhi(f^{*}\circ g)-\varPhi(f)^{*} \circ \varPhi(g).  $$
 Then
 $$ \lvert \phi(Z) \lvert \leq \lvert \phi(X) \lvert ^{\frac{1}{2}} \lvert \phi(Y) \lvert ^{\frac{1}{2}} \,\ \textrm{for all states
 $\phi$ on $\mathbb{A}.$} $$
\end{lemma}
\begin{remark}
  Note that the above inequality holds for Schwarz maps with respect to $C^{*}$ product. It can be shown that completely positive maps with norm less than or equal to 1 are Schwarz maps.
\end{remark}
\textbf{Proof of Theorem \ref{vmoc}:-}
\begin{proof}
 First, we prove that whenever $f,g \in VMO \cap L^{\infty}(\mathbb{T})$, $\{T_{n}(fg)-T_{n}(f)T_{n}(g)\}$ converges to $\{O_{n}\}$ in Type 2 strong cluster sense.\\

Note that for each $n$, $f \mapsto T_{n}(f)$ is a contractive completely positive map. Let $f,g \in VMO \cap L^{\infty}(\mathbb{T})$ and consider the following,\\
$  X_{n}=T_{n}(|f|^{2})-T_{n}(\bar{f})T_{n}(f)\geq 0, \,\ \textrm{(Due to Schwarz inequality)}$\\
 $ Y_{n}=T_{n}(|g|^{2})-T_{n}(\bar{g})T_{n}(g)\geq 0, $\\
 $ Z_{n}=T_{n}(\bar{f}g)-T_{n}(\bar{f})T_{n}(g). $\\
Note that the sequence of matrices $\{X_{n}\},$ $\{Y_{n}\}$ and $\{Z_{n}\}$ are norm bounded.
In particular, we have,
 $ \|Y_{n}\| < \gamma < \infty,$ for all $n \in \mathbb{N}$ and for some $\gamma>0$. By Remark \ref{rem1}, $\{X_{n}\}$ and $\{Y_{n}\}$ converge to $\{O_{n}\}$ in Type 2 strong cluster sense.
 Now for each fixed
$x \in \mathbb{C}^{n}$ with $\|x\|=1$, consider the state $\phi_{x}$ on $\mathcal{B}(\mathbb{C}^{n})$ defined as
$ \phi_{x}(A)=\langle Ax,x\rangle. $
Now by applying Uchiyama's inequality (Lemma \ref{uch}), we get
$$ \lvert \langle Z_{n}x,x \rangle\lvert  \leq \lvert \langle X_{n}x,x\rangle \lvert ^{1/2}  \lvert \langle Y_{n}x,x\rangle \lvert ^{1/2}= \langle X_{n}x,x\rangle^{1/2} \langle Y_{n}x,x\rangle^{1/2} $$
(since $X_{n}$ and $Y_{n}$ are positive, $\langle X_{n}x,x\rangle,\,\  \langle Y_{n}x,x\rangle \geq 0 $ for all $x$).
Let $Z_{n}=B_{n}+iC_{n},$ where $B_{n}$ and $C_{n}$ are self-adjoint for all $n.$ Then,
$$ \lvert \langle B_{n}x,x \rangle\lvert \leq \langle X_{n}x,x\rangle^{1/2} \langle Y_{n}x,x\rangle^{1/2}\,\
\textrm{and}\,\
 \lvert \langle C_{n}x,x \rangle\lvert \leq \langle X_{n}x,x\rangle^{1/2} \langle Y_{n}x,x\rangle^{1/2}, $$
for all $x \in \mathbb{C}^{n}$.
We will show that both $\{B_{n}\}$ and $\{C_{n}\}$ converge to $\{O_{n}\}$ in Type 1 strong cluster sense.

Consider $B_{n}$.
Since $\{X_{n}\}$ converges to $\{O_{n}\}$ in Type 1 strong cluster sense, for $\epsilon>0$, there exist natural numbers $N_{1,\epsilon}$ and $N_{2,\epsilon}$ such that for all $n>N_{2,\epsilon}$, except at most $N_{1,\epsilon}$ eigenvalues, all other eigenvalues are in $[0,\epsilon)$. For any Hermitian matrix $A$, let $\lambda_{k}(A)$ denote the $k^{th}$ eigenvalue of
$A$ when arranged in descending order. Then $\lambda_{k}(X_{n})\in[0,\epsilon)$, for $k > N_{1,\epsilon}$. Also by Min Max principle for eigenvalues of Hermitian 
matrices we have,
$$ \lambda_{k}(X_{n})= \underset{\dim M=k}{\underset{M\subseteq \mathbb{C}^{n}}{\max}}\,\ \underset{x\in M}{\underset{\|x\|=1}{\min}}\langle X_{n}x,x \rangle, $$
$$ [ \lambda_{k}(X_{n})]^{1/2}= \underset{\dim M=k}{\underset{M\subseteq \mathbb{C}^{n}}{\max}}\,\ \underset{x\in M}{\underset{\|x\|=1}{\min}}[\langle X_{n}x,x \rangle]^{1/2}. $$
Now,
 $$ \underset{x \in M}{\underset{\|x\|=1}{\min}} \langle B_{n}x,x \rangle \leq \underset{x \in M}{\underset{\|x\|=1}{\min}} [\langle X_{n}x,x \rangle]^{1/2}\sqrt{\gamma}, $$
for any subspace $M$ of $\mathbb{C}^{n}$. So we obtain,
$$ \underset{\dim M=k}{\underset{M\subseteq H}{\max}}\,\ \underset{x\in M}{\underset{\|x\|=1}{\min}}\langle B_{n}x,x \rangle \leq \sqrt{\epsilon\gamma}. $$
Let $\delta>0$ and put $\epsilon=\frac{\delta^{2}}{\gamma}$, then $\lambda_{k}(B_{n}) \leq \delta $ for $k> N_{1,\epsilon}$.
Similarly, $\lambda_{k}(-B_{n})\leq \delta$, for $k>N_{1,\epsilon}$.
But $\lambda_{k}(-B_{n})=-\lambda_{n-k+1}(B_{n})$. Therefore $\lambda_{k}(B_{n}) \leq \delta $ for $k= N_{1,\epsilon}+1,\hdots,n$ and $\lambda_{k}(B_{n}) \geq - \delta $
except at most for $k=n-N_{1,\epsilon}+1,\hdots,n$. This will imply that $-\delta \leq \lambda_{k}(B_{n}) \leq \delta $, for $k=N_{1,\epsilon}+1,\hdots,n-N_{1,\epsilon}$. Hence except
at most $2N_{1,\epsilon}$ eigenvalues, all other eigenvalues of $B_{n}$ are in $[-\delta,\delta]$ and note that $N_{1,\epsilon}=O(1)$. Thus $\{B_{n}\}$ converges to
$\{O_{n}\}$ in Type 1 strong cluster sense. Similarly $\{C_{n}\}$ converges to $\{O_{n}\}$ in Type 1 strong cluster sense. Since both are self-adjoint, they converge to $\{O_{n}\}$ in Type 2 strong cluster sense also. Hence, $\{Z_{n}\}$ converges to $\{O_{n}\}$ in Type 2
strong cluster sense. Replacing $f$ by $\bar{f}$, we obtain that $\{T_{n}(fg)-T_{n}(f)T_{n}(g)\}$ converges to $\{O_{n}\}$ in Type 2 strong cluster sense.\\
 
 Let $M$ be a $C^{*}$-subalgebra of $L^{\infty}(\mathbb{T})$ such that whenever $f,g \in M$, $\{T_{n}(fg)-T_{n}(f)T_{n}(g)\}$ converges to $\{O_{n}\}$ in Type 2 strong cluster sense. Let $f \in M$. Then $\{T_{n}(|f|^{2})-T_{n}(f)T_{n}(\bar{f})\}$ and $\{T_{n}(|f|^{2})-T_{n}(\bar{f})T_{n}(f)\}$ converge to $\{O_{n}\}$ in Type 2 strong cluster sense. Hence by Theorem \ref{vmo}, $f \in VMO \cap L^{\infty}(\mathbb{T})$. Thus $M \subseteq VMO \cap L^{\infty}(\mathbb{T})$.
\end{proof}
\begin{remark}
 To prove $M \subseteq VMO \cap L^{\infty}(\mathbb{T})$, we have only used the fact that whenever $f \in M$, $\{T_{n}(|f|^{2})-T_{n}(f)T_{n}(\bar{f})\}$ and $\{T_{n}(|f|^{2})-T_{n}(\bar{f})T_{n}(f)\}$ converge to $\{O_{n}\}$ in Type 2 strong cluster sense. By similar arguements, we obtain the following. Let $M_{1}$ be the set of all $f \in L^{\infty}(\mathbb{T})$ such that for $f \in M_{1}$, $\{T_{n}(|f|^{2})-T_{n}(f)T_{n}(\bar{f})\}$ and $\{T_{n}(|f|^{2})-T_{n}(\bar{f})T_{n}(f)\}$ converge to $\{O_{n}\}$ in Type 2 strong cluster sense. Then $M_{1}=VMO \cap L^{\infty}(\mathbb{T})$. 
\end{remark}
\section{Concluding remarks and future problems}
Here we could achieve the following.
\begin{enumerate}
 \item As mentioned, the problem of checking the compactness of $T_{|f|^{2}}-T_{f}T_{\bar{f}}$ is one of the important problems in operator theory. We could characterize the compactness of $T_{|f|^{2}}-T_{f}T_{\bar{f}}$ and $T_{|\tilde{f}|^{2}}-T_{\tilde{f}}T_{\bar{\tilde{f}}}$ with the convergence of $\{T_{n}(|f|^{2})-T_{n}(f)T_{n}(\bar{f})\}$ in the sense of singular value clustering (for $f \in BMO$). The advantage is that we can compute the matrix $T_{n}(|f|^{2})-T_{n}(f)T_{n}(\bar{f})$ using the Fourier coefficients of $f$.
 \item The largest $C^{*}$ subalgebra $M$ of $L^{\infty}(T)$ such that whenever $f,g \in M$, $\{T_{n}(fg)-T_{n}(f)T_{n}(g)\}$ converges to $\{O_{n}\}$ in Type 2 strong cluster sense coincides with $VMO \cap L^{\infty}(\mathbb{T})$, the largest $C^{*}$-subalgebra of $L^{\infty}(\mathbb{T})$ such that whenever $f,g \in VMO \cap L^{\infty}(\mathbb{T})$, $T_{fg}-T_{f}T_{g}$ is compact.
 \item As an application of Theorem \ref{vmoc}, we can extend the Korovkin-type theorem (Theorem 3.4 of \cite{laa}) obtained for Toeplitz operators with continuous symbols to Toeplitz operators with $VMO \cap L^{\infty}(\mathbb{T})$ symbols.
 \item It is known that for $f,g \in L^{2}(\mathbb{T})$, $\{T_{n}(fg)-T_{n}(f)T_{n}(g)\}$ converges to $\{O_{n}\}$ in Type 2 weak cluster sense (see Theorem 2.9 and 7.8 of \cite{serb}). We have obtained that the convergence is in strong cluster sense if $f \in VMO$ and $g=\bar{f}$. One shall check whether the class of functions $f$ for which $\{T_{n}(fg)-T_{n}(f)T_{n}(\bar{f})\}$ converges to $\{O_{n}\}$ in Type 2 strong cluster sense can be characterized by $VMO$.
 \item We shall look for results analogous to Theorem \ref{eq}, Theorem \ref{vmo} and Theorem \ref{vmoc} in the case of Toeplitz operators acting on Bergman space, Fock space etc.
\end{enumerate}
\textbf{Acknowledgement :-} The author wishes to thank Wolfram Bauer and V. B. Kiran Kumar  for fruitful discussions and valuable suggestions. The author is supported by the Post Doctoral Fellowship scheme of Cochin University of Science and Technology.\\

\nocite{*}
\bibliographystyle{amsplain}
\bibliography{rahul_arxiv}

\end{document}